\documentclass[11pt]{article}
\usepackage{amsmath,amsfonts,amssymb,theorem,graphics,verbatim}

\def\BEN{\begin{enumerate}}  \def\BI{\begin{itemize}}
\def\EEN{\end{enumerate}}   \def\EI{\end{itemize}}
    \def\nn{\nonumber}

\def\beq{\begin{eqnarray}} \def\eeq{\end{eqnarray}}

\def\eqn#1{\begin{equation}#1\end{equation}}

\def\al*#1{\begin{align*}#1\end{align*}}

\def\ga*#1{\begin{gather*}#1\end{gather*}}

\def\alat*#1#2{\begin{alignat*}{#1}#2\end{alignat*}}
\def\bea{\begin{eqnarray*}}
\def\eea{\end{eqnarray*}}
\def\ml*#1{\begin{multline*}#1\end{multline*}}

 \def\mbf{\mathbf}

 \def\le{\left} \def\ri{\right} \def\i{\infty}
   \def\R{{\mathbb R}}

\def\FF{\mathcal{F}}

\def\te#1{\mathrm{e}^{#1}} \def\T{\tilde}  

\def\WH{\widehat}

     \def\a{\alpha} \def\b{\beta}
     \def\th{\theta}

 \def\k{\kappa}   
  \def\nn{\nonumber}   \def\s{\sigma}
\def\t{\tau}     
  \def\q{\qquad} 
 \def\G{\Gamma}   
  \def\td{\text{\rm d}}
\numberwithin{equation}{section}

\newtheorem{Thm}{Theorem}
 \newtheorem{Lemma}{Lemma}
\newtheorem{Prop}{Proposition} 
{\theorembodyfont{\normalfont}
 
\newtheorem{Rem}{Remark} 
 
}

\newcommand{\exit}{{\mbox{\, \vspace{3mm}}} \hfill\mbox{$\square$}}

\def\sh{\Gamma} 
\begin{document}
\title{Cram\'{e}r asymptotics for finite time first passage
probabilities of general L\'{e}vy processes}
\author{Zbigniew Palmowski\footnote{University of Wroc\l aw, pl. Grunwaldzki
2/4, 50-384 Wroc\l aw, Poland, E-mail: zpalma@gmail.com} \q
Martijn Pistorius\footnote{Imperial College London,
Department of Mathematics, South Kensington Campus, London SW7
2AZ, UK, E-mail:m.pistorius@imperial.ac.uk}} \maketitle
 \begin{abstract}
We derive the exact asymptotics of  $P(\sup_{u\leq t}X(u)
> x)$ if $x$ and $t$ tend
to infinity with $x/t$ constant, for a general L\'{e}vy process $X$ that
admits exponential moments. The proof is based on a renewal
argument and a two-dimensional renewal theorem of H\"{o}glund \cite{H}.
\end{abstract}

\section{Introduction}

The study of boundary crossing probabilities of L\'{e}vy processes
has applications in many fields, including ruin theory (see e.g.
Rolski et al. \cite{RSST} and Asmussen \cite{Asm2000}), queueing
theory (see e.g. Borovkov \cite{Borovkov} and Prabhu
\cite{Prabhu}), statistics (see e.g. Siegmund \cite{siegmund}) and
mathematical finance (see e.g. Roberts and Shortland \cite{RS}).

As in many cases closed form expressions
for (finite time) first passage probabilities
are either not available or intractable, a good deal of the literature
has been devoted to logarithmic or exact asymptotics for
first passage probabilities, using different techniques. Martin-L\"{o}f \cite{mar}
and Collamore \cite{C} derived large deviation
results for first passage probabilities of a general class of processes.
Employing two-dimensional renewal theory and asymptotic
properties of ladder processes, respectively, H\"{o}glund \cite{H} and
von Bahr \cite{bahr} obtained exact asymptotics for ruin probabilities of
the classical risk process (see also Asmussen \cite{Asm2000}).
Bertoin and Doney \cite{BD} generalised the classical
Cram\'{e}r-Lundberg approximation (of the perpetual ruin probability
of a classical risk process) to general L\'{e}vy processes.

In this paper we obtain the exact asymptotics of the finite time ruin probability
$P(\tau(x)\leq t)$, where $\tau(x)=\inf\{t\geq 0: X(t)>x\}$, for a
general L\'{e}vy process $X(t)$ ($X(0)=0$), if $x$ and
$t$ jointly tend to infinity in fixed proportion,
generalising Arfwedson \cite{Arfwedson} and H\"{o}glund \cite{H}
who treated the case of a classical risk process.
The proof is based on an embedding of the ladder process of $X$
and a two-dimensional renewal theorem of H\"{o}glund \cite{H}.

The remainder of the paper is organized as follows. In Section
\ref{mainresult} the main result is presented, and its proof
is given in Section \ref{sec:thm}.

\section{Main result}\label{mainresult}

Let $X$ be a L\'{e}vy process with non-monotone paths that
satisfies
\begin{equation}\label{eq:nu}
E[\te{\alpha_0 X(1)}]<\i \quad \text{for some $\alpha_0>0$},
\end{equation}
and denote by $\t(x) = \inf\{t\ge0: X(t) > x\}$ the first crossing
time of $x$. We exclude the case that $X$ is a compound Poisson
process with non-positive infinitesimal drift, as this corresponds
to the random walk case which has already been treated in the
literature.

The law of $X$ is determined by its Laplace exponent $\psi(\th) =
\log E[\te{\th X(1)}]$ that is well defined on the maximal domain
$\Theta = \{\theta\in\R: \psi(\theta)<\i\}$. Restricted to the
interior $\Theta^o$, the map $\th\mapsto \psi(\th)$ is convex and
differentiable, with derivative $\psi'(\th)$.\footnote{For
$\th\in\Theta\backslash\Theta^o$, $\psi'(\theta)$ is understood to
be $\lim_{\eta\to\theta,\eta\in\Theta^o}\psi'(\eta)$. } Moreover,
$\psi'(0+)=E[X(1)]$ if $E[|X(1)|]<\infty$.
By the strict convexity of $\psi$, it follows that $\psi'$ is
strictly increasing on $(0,\i)$ and we denote by
$\sh:\psi'(0,\i)\to (0,\i)$ its right-inverse function.

Associated to the measure $P$ is the exponential family of
measures $\{P^{(c)}: \text{$c\in\Theta$}\}$ defined by their
Radon-Nikodym derivatives
\eqn{\label{eq:changeofmeasure} \le.\frac{\td P^{(c)}}{\td
P}\ri|_{\FF_t} = \exp\le(c X(t) - \psi(c)t\ri). }
It is well known that under this change of measure X is still a
L\'{e}vy processes and its new Laplace exponent satisfies
\begin{equation}\label{newkappa}\psi^{(c)}(\alpha)=\psi(\alpha+c)-\psi(c).\end{equation}

Related to $X$ and its running supremum are the local time $L$ of
$X$ at its supremum, its right-continuous inverse $L^{-1}$ and the
upcrossing ladder process $H$ respectively. The Laplace exponent
$\kappa$ of the bivariate (possibly killed) subordinator
$(L^{-1},H)$,
\begin{equation}\label{kappa}
\te{-\k(\a,\b)t} = E[\te{-\a L^{-1}_t -\b H_t}\mbf
1_{(L^{-1}_t<\i)}],
\end{equation}
is related to $\psi$ via the Wiener-Hopf factorisation identity
\begin{equation}\label{wh}
u-\psi(\theta) =  k \kappa(u,- \theta)\WH\kappa(u, \theta), \quad u\ge 0, \theta\in\Theta^o,
\end{equation}
for some constant $k>0$ where $\WH\kappa$ is the Laplace exponent
of the dual ladder process. Refer to Bertoin \cite[Ch. VI]{B} for
further background on the fluctuation theory of L\'{e}vy
processes.

Bertoin and Doney \cite{BD} showed that, if the Cram\'{e}r condition
holds, that is $\gamma>0$, where \begin{equation}\gamma :=
\sup\{\th\in\Theta:\psi(\th) = 0\},\end{equation} the
Cram\'{e}r-Lundberg approximation remains valid for a general
L\'{e}vy process:
\begin{equation}\label{c}
\lim_{x\to\i}\te{\gamma x}P[\tau(x) < \i] = C_\gamma,
\end{equation}
where $C_\gamma\ge 0$ is positive if and only if $E[\te{\gamma
X(1)}|X(1)|]<\i$ and is then given by $C_\gamma =
\beta_\gamma/[\gamma m_\gamma]$, where
$$\beta_\gamma = -\log P[H_1 < \i], \quad
m_\gamma = E[\te{\gamma H_1}H_1\mbf 1_{(H_1<\i)}].$$ Further,
Doob's optional stopping theorem implies the following bound:
\begin{equation}\label{eq:c1}
\te{\gamma x}P(\tau(x)<\infty)=
E^{(\gamma)}[\te{-\gamma(X(\tau(x))-x)}\mbf
1_{(\tau(x)<\infty)}]\leq 1.
\end{equation}

The result below concerns the asymptotics of the finite time ruin
probability $P(\tau(x)\leq t)$ when $x,t$ jointly tend to infinity
in fixed proportion. For a given proportion $v$ the rate of decay
is either equal to $\gamma v t$ or to $\psi^*(v)t$, where $\psi^*$
is the convex conjugate of $\psi$:
$$\psi^*(u) = \sup_{\alpha\in\R}(\alpha u - \psi(\alpha)).
$$
We restrict ourselves to L\'{e}vy processes satisfying the
following condition
\begin{equation}
\tag*{(H)} \text{$\sigma>0$ or the L\'{e}vy measure is
non-lattice,}
\end{equation}
where $\sigma$ denotes the Gaussian coefficient of $X$. Recall
that a measure is called non-lattice if its support is not
contained in a set of the form $\{a + bh, h\in\mathbb Z\}$, for
some $a,b>0$. Note that (H) is satisfied by any L\'{e}vy process
whose L\'{e}vy measure has infinite mass.

We write $f\sim g$ if $\lim_{x,t\to\infty, x=vt+{\rm o}(t^{1/2})} f(x,t)/g(x,t) =
1$.

\begin{Thm}\label{thm}
Assume that $(H)$ holds. Suppose that $0<\psi'(\gamma)<\i$ and
that there exists a $\sh(v)\in\Theta^\circ$ such that
$\psi'(\sh(v))=v$. If $x$ and $t$ tend to infinity such that $x=
vt+{\rm o}(t^{1/2})$ then
$$
P(\t(x) \leq t) \sim \begin{cases} C_\gamma\te{-\gamma x},
& \text{if $0< v < \psi'(\gamma)$,}\\
D_v t^{-1/2}\te{-\psi^*(v)t}, & \text{if $ v > \psi'(\gamma)$,}
\end{cases}
$$
with $C_0=1$ and $D_v$ given by
\begin{eqnarray*}
D_v &=& \frac{- v \log E[\te{-\eta_v L_1^{-1}}\mbf 1_{(L_1^{-1}<\i)}]}%
{\eta_v E[\te{\sh(v) H_1-\eta_v L_1^{-1}}H_1\mbf
1_{(L_1^{-1}<\i)}]} \times \frac{1}{\sh(v)
\sqrt{2\pi\psi''(\sh(v))}},
\end{eqnarray*}
where $\eta_v = \psi(\sh(v))$.
\end{Thm}

\begin{Rem}\label{uwaga1}
(a) For a spectrally negative L\'{e}vy process the joint exponent
of the ladder process is given by $\kappa(\alpha,
\beta)=\beta+\Phi(\alpha)$ ($\alpha, \beta\geq 0$), where
$\Phi(\a)$ is the largest root of $\psi(\th)=\a$, and thus
\begin{equation}\label{cv}D_v=\overline{D}_v:=\frac{v}{ \psi(\sh(v))\sqrt{2\pi\psi''(\sh(v))}}, \quad
C_\gamma \equiv 1.
\end{equation}
Indeed, \begin{eqnarray*}&&D_v=\overline{D}_v\times
\frac{\kappa(\eta_v,0)}{\Gamma(v)\frac{\partial}{\partial
\beta}\kappa(\eta_v,
\beta)_{|\beta=-\Gamma(v)}\exp\{-\kappa(\eta_v,-\Gamma(v))\}}\\
&&\ \ \ \ =\overline{D}_v\times
\frac{1}{\exp\{-\Phi(\eta_v)+\Gamma(v)\}}=\overline{D}_v\end{eqnarray*}
since $\Phi(\eta_v)=\Gamma(v)$.

(b) If $X$ is spectrally positive, $\kappa(\alpha,\beta) = [\alpha
- \psi(-\beta)]/[\hat\Phi(\alpha) - \beta]$ (see e.g. \cite[Thm
VII.4]{B}), where $\hat\Phi(\a)$ is the largest root of
$\psi(-\th)=\a$ and we find that
$$
D_v = \frac{\sh(v) + \T\sh(v)}{\sh(v)\T\sh(v)}
\frac{1}{\sqrt{2\pi\psi''(\sh(v))}}, \quad C_\gamma =
\frac{\psi'(0)}{\psi'(\gamma)},
$$
where $\T\sh(v) = \sup\{\th:\psi(-\th) = \psi(\G(v))\}$,
recovering formulas that can be found in Arfwedson \cite{Arfwedson} and
Feller \cite{Feller} respectively, for the case of a classical risk
process.
\end{Rem}

\begin{Rem}\label{uwaga2}
Heuristically, in the case $v>\psi'(\gamma)$, the asymptotics in
Thm. \ref{thm} can be regarded as a consequence of the central limit
theorem, that is, under the tilted measure $P^{\Gamma(v)}$,
asymptotically
$$\frac{\tau(x)-x/v}{\omega \sqrt{x}}$$
 follows a standard normal distribution, where by
(\ref{newkappa}) and choice of $\sh(v)$,
$$\omega^2=\frac{{\rm
Var}^{(\sh(v))}[X_1]}{\left(E^{(\sh(v))}[X_1]\right)^3}=\frac{\psi^{(\sh(v))\prime\prime}(0)
}{\left(\psi^{(\sh(v))\prime}(0)\right)^3}=
\frac{\psi''(\sh(v))}{v^3}.
$$
This explains why the asymptotics remain valid if $x$ deviates
${\rm o}(x^{1/2})={\rm o}(t^{1/2})$ from the line $vt$.

In the boundary case $v=\psi'(\gamma)$, in which case
$E^{(\Gamma(v))}[\tau(x)]=t$, the exact asymptotics of
$P(\t(x) \leq t)$ may depend on the way in which $x/t$ tends to $v$.
Note that this case is excluded from Theorem \ref{thm}.
\end{Rem}

\begin{Rem}\label{rem:sp}
In the case $0<v<\psi'(\gamma)$, the asymptotics in Theorem
\ref{thm} are a consequence of the law of large numbers. To see
why this is the case, note that $\te{\gamma x}P(\tau(x)\leq t) =
\te{\gamma x}P(\tau(x) < \i) - \te{\gamma x}P(t<\tau(x) < \i)$,
where the first term tends to $C_\gamma$ in view of (\ref{c}),
while for the second term the Markov property and \eqref{eq:c1}
imply that
\begin{eqnarray*}
\lefteqn{\te{\gamma x}P(t<\tau(x) < \i)}\\&&=\int_{-\infty}^x
P(\tau(x)>t, X(t)\in dy)\te{\gamma y}\te{\gamma (x-y)}P(\tau(x-y)<\infty)\\
&&\leq \int_{-\infty}^x P(X(t)\in dy)\te{\gamma
y}=P^{(\gamma)}(X(t)\leq x),
\end{eqnarray*}
which tends to $0$ as $t$ tends to infinity in view of the law of
large numbers since $E^{(\gamma)}[X(t)]=t\psi'(\gamma)>x$. The proof below deals with the case that $v>\psi'(\gamma)$.
\end{Rem}

\section{Proof of Theorem \ref{thm}}\label{sec:thm}
The idea of the proof is to lift asymptotic results that have been
established for random walks by H\"oglund \cite{H} and Arfwedson
\cite{Arfwedson} to the setting of L\'{e}vy processes by considering
suitable random walks embedded in the L\'{e}vy process (more
precisely, in its ladder process). We first briefly recall these
results following the H\"{o}glund \cite{H} formulation.

\subsection{Review of H\"oglund's random walk asymptotics}
Let $(S,R)= \{(S_i,R_i), i=1,2,\ldots\}$ be a (possibly killed)
random walk starting from $(0,0)$ whose components $S$ and
$R$ have non-negative increments, and consider the crossing
probabilities
\begin{eqnarray*}
G_{a,b}(x,y) &=& P(N(x)<\i, S_{N(x)} > x + a, R_{N(x)} \leq y+b),  \\
K_{a,b}(x,y) &=& P(N(x)<\i, S_{N(x)} > x + a, R_{N(x)} \geq y+b),
\end{eqnarray*}
where $a\ge 0, b\in\R$ and $N(x) = \min\{n: S_n  > x\}$. Let $F$
denote the (possibly defective) distribution function of the
increments of the random walk with joint Laplace transform $\phi$
and set $F_{(u,v)}(\td x, \td y) = \te{-ux - vy}F(\td x,\td
y)/\phi(u,v)$. Let $$V(\zeta)=E_\zeta[(R_1E_\zeta[S_1] - S_1
E_\zeta[R_1])^2]/ E_{\zeta}[S_1]^3$$ for $\zeta=(\xi,\eta)$ where
$E_\zeta$ denotes the expectation w.r.t. $F_\zeta$.

For our purposes it will suffice to consider
random walks that satisfy the following
non-lattice assumption (the analogue of the non-lattice
assumption in one dimension):
\begin{equation*}
\tag*{(G)}
\text{The additive group spanned by the support of $F$ contains $\R^2_+$.}
\end{equation*}
Specialised to our setting Prop. 3.2 in H\"{o}glund (1990) jointly
with the remark given on p. 380 therein read as follows:
\begin{Prop}\label{hog}
Assume that $(G)$ holds, and that there exists a
$\zeta=(\xi,\eta)$ with $\phi(\zeta)=1$ such that $v =
E_\zeta[S_1]/E_\zeta[R_1]$, where $\phi$ is finite in a
neighbourhood of $\zeta$ and $(0,\eta)$. If $x,y$ tend to infinity
such that $x=vy+{\rm o}(y^{1/2})>0$ then it holds that
\begin{eqnarray*}
G_{a,b}(x,y) &\sim& D(a,b) x^{-1/2}\te{x\xi+y\eta}
\q \text{if\, $\eta>0$},\\
K_{a,b}(x,y) &\sim& D(a,b) x^{-1/2}\te{x\xi+y\eta} \q \text{if\,
$\eta<0$},
\end{eqnarray*}
for $a\ge 0, b\in\R$, where $D(a,b) = C(a,b) \cdot (2\pi
V(\zeta))^{-1/2} $, with $V(\zeta)>0$ and
$$
C(a,b) = \frac{1}{|\eta| E_\zeta[S_1]} \te{b\eta}\int_a^\i
P_\zeta(S_1 \ge x)\te{\xi x}\td x.
$$
\end{Prop}

\subsection{Embedded random walk}
Denote by $e_1, e_2, \ldots$ a sequence of independent exp$(q)$
distributed random variables and by $\s_n=\sum_{i=1}^n e_i$, with
$\s_0=0$, the corresponding partial sums, and consider the
two-dimensional (killed) random walk $\{(S_i,R_i), i=1,2\ldots\}$
starting from $(0,0)$ with step-sizes distributed according to
$$
F^{(q)}(\td t, \td x) = P(H_{\s_{1}}\in \td x,
L^{-1}_{\s_{1}}\in\td t),
$$
and write $G^{(q)}$ for the corresponding crossing probability
$$
G^{(q)}(x,y) = G_{0,0}(x,y) = P(N(x)<\i, R_{N(x)}
\leq y).
$$
Note that $F^{(q)}$ is a probability measure that is defective
precisely if $X$ drifts to $-\i$, with Laplace transform $\phi$
given by
$$
\phi(u,v) = \iint \te{-ut-vx} F^{(q)}(\td t, \td x)  = \frac{q}{q
- \k(u,v)}.
$$

The key step in the proof is to derive bounds for $P(\tau(x)\leq
t)$ in terms of crossing probabilities involving the random walk
$(S,R)$:

\begin{Lemma}\label{lem:est} Let $M,q>0$. For $x,t>0$ it holds that
\begin{equation}\label{star}
G^{(q)}(x,t) \leq P(\tau(x)\leq t) \leq G^{(q)}(x,t+M)/h(0-, M),
\end{equation}
where $h(0-,M) = \lim_{x\uparrow 0} h(x,M)$, with $h(x,t) :=
P(H_{\s_1} > x, L^{-1}_{\s_1} \leq t)$.
\end{Lemma}

{\it Proof:} Let $T(x) = \inf\{t\ge0: H_t > x\}$ and note that
$\t(x) = L^{-1}_{T(x)}$. By applying the Markov property it
follows that
\begin{eqnarray}
\nn  P(\t(x)\leq t) &=& P(T(x) < \i, L^{-1}_{T(x)} \leq t)\\
  &=& \sum_{n=1}^\i P(\s_{n-1} \leq T(x) < \sigma_n, L^{-1}_{T(x)} \leq
  t)\\
\nn  &=& \sum_{n=1}^\i P(H_{\s_{n-1}}\leq x, H_{\s_n} > x,
L^{-1}_{T(x)} \leq
  t)\\
\nn  &=& \sum_{n=1}^\i\iint  P(H_{\s_{n-1}}\in \td y,
  L^{-1}_{\s_{n-1}}\in\td s) \\&& \quad \quad \times\ P(H_{\s_1} > x - y, L^{-1}_{T(x-y)} \leq
  t-s)\\
  &=& \sum_{n=0}^\i F^{(q)\star n} \star f(x,t) = (U \star
  f)(x,t),\label{eq:u*f}
\end{eqnarray}
where $U=\sum_{n=0}^\infty F^{(q)\star n}$, $f(x,t) = P(H_{\s_1}
> x, L^{-1}_{T(x)} \leq t)$ and $\star$ denotes convolution.
Following a similar reasoning it can be
checked that \begin{equation}\label{eq:u*h}G^{(q)}(x,t) = U\star
h(x,t).\end{equation} In view of \eqref{eq:u*f} and
\eqref{eq:u*h}, the lower bound in \eqref{star} follows since $$f(x,t) \ge h(x,t),$$ taking note of the fact that $H_{\s_1}> x$
precisely if $T(x) < \s_1$,
while the upper bound in \eqref{star} follows by observing that
for fixed $M>0$,
\begin{eqnarray*}
\nn  h(x,t+M) &\ge& P(H_{\s_1} > x, L^{-1}_{T(x)}
\leq t, L^{-1}_{\s_1} - L^{-1}_{T(x)} \leq M)\\
\nn &=& P(H_{\s_1} > x, L^{-1}_{T(x)} \leq t) P(L^{-1}_{\s_1} \leq
M)\\ &=& f(x,t)h(0-,M),
\end{eqnarray*}
where we used the strong Markov property of $L^{-1}$ and the lack
of memory property of $\s_1$. \exit

\medskip

\noindent Applying H\"oglund's asymptotics in Proposition \ref{hog} yields
the following result:

\begin{Lemma}\label{hog2}
Let the assumptions of Proposition \ref{hog} hold true. If
$x,t\to\i$ such that for $v > \psi'(\gamma)$ we have $x=vt+{\rm o}(t^{1/2})$ then
$$
G^{(q)}(x,t+M) \sim D_{q,M}t^{-1/2}\te{-\psi^*(v)t}, \quad M\ge 0,
$$
where $D_{q,M} = \frac{v}{\sqrt{2\pi\psi''(\sh(v))}} C_{q,M}$ with
$$
C_{q,M} =
\te{\psi(\sh(v))M}\frac{\k(\psi(\sh(v)),0)}{c_v\psi(\sh(v))\sh(v)}
\frac{q}{q+\k(\psi(\sh(v)),0)},
$$
where $c_v=E[\te{\sh(v) H_{1} - \psi(\G(v))L^{-1}_{1}}H_{1} \mbf
1_{(L^{-1}_{1}<\i)}]$.
\end{Lemma}
Lemma \ref{hog2} is a consequence of the following auxiliary
identities:

\begin{Lemma}\label{lem:mom} Let $u>\gamma$, $u\in\Theta^o$.
\begin{eqnarray}
\label{eq:a}\phi(z,-u)&=&1 \q\text{iff  $\k(z,-u)=0$ iff
$\psi(u)=z$}
\\
\label{eq:b} \psi'(u) &=&
E^{(u)}[X(1)]=E^{(u)}[H_{\s_1}]\cdot(E^{(u)}[L^{-1}_{\s_1}])^{-1}\\
\nn  \psi''(u)& =& E^{(u)}[(H_{\s_1} -
\psi'(u)L^{-1}_{\s_1})^2]\cdot(E^{(u)}[L^{-1}_{\s_1}])^{-1}\\
\label{eq:c} &=&\psi'(u)E^{(u)}[(H_{\s_1} -
\psi'(u)L^{-1}_{\s_1})^2]\cdot(E^{(u)}[H_{\s_1}])^{-1}\\
\label{eq:d} \psi^*(v) &=& v\Gamma(v) - \psi(\Gamma(v))\q
\text{for $v>0$ with $\Gamma(v)\in\Theta^o$.}
\end{eqnarray}

\end{Lemma}

{\it Proof}: Eq (\ref{eq:a}): Note that for $u,z>0$ it holds that
 $\WH\kappa(z,u) > 0$. In view of the identity \eqref{wh}
the statement follows.

Eq (\ref{eq:b}): Note that if $u>\gamma$ then by the fact that
$\psi(0)=\psi(\gamma) = 0$ and the strict convexity of $\psi$ it
follows that $\psi(u) > 0$. In view of \eqref{wh} it follows then
that $\k(\psi(u),-u)=0$ for $u\in\Theta^o$, $u>\gamma$.
Differentiating with respect to $u$ shows that
\begin{equation}\label{eq:psip}
\psi'(u) =
\partial_2\k(\psi(u),-u)(\partial_1\k(\psi(u),-u))^{-1}.
\end{equation}
Also, note that $E^{(u)}[H_{\s_1}] = q^{-1}E^{(u)}[H_1]$,
$E^{(u)}[L^{-1}_{\s_1}] = q^{-1}E^{(u)}[L^{-1}_1]$ and
\begin{equation*}
E^{(u)}[H_1]=\partial_2\k(\psi(u),-u), \quad
E^{(u)}[L^{-1}_1]=\partial_1\k(\psi(u),-u).
\end{equation*}
Eq (\ref{eq:c}) follows as a matter of calculus, by
differentiation of \eqref{eq:psip} with respect to $u$. Finally,
Eq. \eqref{eq:d} follows from the definition of $\psi^*$. \exit

{\it Proof of Lemma \ref{hog2}} The proof follows by an
application of Prop. \ref{hog} to $G^{(q)}(x,t+M)$ with
$$(S_1, R_1)
= (H_{\s_1},L^{-1}_{\s_1})\q \text{ and }\q
\zeta=(-\Gamma(v),\eta_v).$$ Note that, by \eqref{eq:a} with
$u=\Gamma(v)$, $\phi(\zeta)=1$, and that $\eta_v = \psi(\G(v))>0$
if $v > \psi'(\gamma)$. For this choice of the parameters,
$E_\zeta [S_1]=E^{(\Gamma(v))}[H_{\s_1}]=c_v/q$, and Eqs.
\eqref{eq:d}, \eqref{eq:b},\eqref{eq:c} imply that  $\xi x+\eta
t=-\psi^*(v)t$ and
$$V(\zeta)=
\psi''(\Gamma(v))/\psi'(\Gamma(v))=\psi''(\Gamma(v))/v.$$ To
complete the proof we are left to verify
 the form of the constants. The calculation of the
$C_{q,M}=C(0,0)\te{\eta M}$ goes as follows:
\begin{eqnarray*}
C_{q,M} &=&
\frac{q\te{\psi(\sh(v))M}}{\psi(\sh(v))c_v}\le(\int_0^\i
\te{-\sh(v)x}E[\te{\sh(v) H_{\s_1} - \psi(\G(v))L^{-1}_{\s_1}}
\mbf 1_{(x\leq H_{\s_1} < \i)}]\td x \ri)\\
&=& \frac{q\te{\psi(\sh(v))M}}{\psi(\sh(v))\sh(v)c_v}
\le(1 - E[\te{-\psi(\sh(v))L^{-1}_{\s_1}}\mbf 1_{(L^{-1}_{\s_1} < \i)}]\ri)\\
&=& \frac{q\te{\psi(\sh(v))M}}{\psi(\sh(v))\sh(v)c_v}\le(1 -
\frac{q}{q+\k(\psi(\sh(v)),0)}\ri)\\
&=& \frac{q\te{\psi(\sh(v))M}}{\psi(\sh(v))\sh(v)c_v}
\frac{\k(\psi(\sh(v)),0)}{q+\k(\psi(\sh(v)),0)},
\end{eqnarray*}
in view of the definition \eqref{kappa} of $\kappa$. Combining all
results completes the proof.\exit

As final preparation for the proof of Theorem \ref{thm}
we show that the non-lattice condition holds:

\begin{Lemma}
Suppose that $(H)$ holds true. Then $F^{(q)}$ satisfies (G).
\end{Lemma}

{\it Proof}: The assertion is a consequence of the following
identity between measures on $(0,\i)^2$ (which is itself a
consequence of the Wiener-Hopf factorisation, see e.g. Bertoin
\cite[Cor VI.10]{B})
\begin{equation}\label{eq:wh}
P(X_t\in\td x)\td t = t\int_0^\i P(L^{-1}_u\in\td t, H_u\in\td x)
u^{-1}\td u.
\end{equation}
Fix $(y,v)\in (0,\i)^2$ in the support of $\mu_X(\td t,\td x) =
P(X_t\in\td x)\td t$ and let $B$ be an arbitrary open ball around
$(y,v)$. Then $\mu_X(B)>0$; in view of the identity \eqref{eq:wh}
it follows that there exists a set $A$ with positive Lebesgue
measure such that $P((L^{-1}_u, H_u)\in B) > 0$ for all $u\in A$
and thus $P((L^{-1}_{\s_1}, H_{\s_1})\in B) > 0$. Since $B$ was
arbitrary we conclude that $(y,v)$ lies in the support of
$F^{(q)}$. To complete the proof we next verify that if a L\'{e}vy
process $X$ satisfies (H) then $\mu_X$ satisfies (G). To this end,
let $X$ satisfy (H). Suppose first that its L\'{e}vy measure $\nu$
has infinite mass or $\sigma>0$. Then $P(X_t=x)=0$ for any $t>0$
and $x\in\R$, according to Sato \cite[Thm. 27.4 ]{Sato}. Thus, the
support of $P(X_t\in\td x)$ is uncountable for any $t>0$, so that
$\mu_X$ satisfies (G). If $\nu$ has finite mass then it is
straightforward to verify that $P(X_t\in\td x)$ is non-lattice for
any $t>0$ if $\nu$ is, and that then $\mu_X$ satisfies (G). \exit
\bigskip

{\it Proof of Theorem \ref{thm}:} Suppose that $v > \psi'(\gamma)$
(the case $v<\psi'(\gamma)$ was shown in Remark \ref{rem:sp}).
Writing $l(t,x) = t^{1/2}\te{\psi^*(v)t}P(\tau(x)\leq t)$, Lemmas
\ref{lem:est}, \ref{hog2}  and \ref{lem:mom} imply that
\begin{eqnarray*}
s &=& \limsup_{x,t\to\i, x=tv+{\rm o}(t^{1/2})} l(t,x) \leq D_{q,M}/h(0-,M),\\
i &=& \liminf_{x,t\to\i, x=tv+{\rm o}(t^{1/2})} l(t,x) \geq D_{q,0}.
\end{eqnarray*}
By definition of $h$ and $D_{q,M}$ it directly follows that, as
$q\to\i$,
$$D_{q,0} \to D_v,\ D_{q,M} \to D_v \te{\psi(\G(v))M}\text{ and }
h(0-,M)=P(L^{-1}_{\sigma_1}\leq M)\to 1.$$ Letting $M\downarrow 0$
yields that $s=i = D_v$, and the proof is complete. \exit

\subsection*{Acknowledgments}
We thank Florin Avram for helpful advice and useful comments. We
are grateful to the referee for supplying the short direct proof
(Remark \ref{rem:sp}) for the case $v<\psi'(\gamma)$. This work is
partially supported by EPSRC grant EP/D039053/1 and the Ministry
of Science and Higher Education of Poland under the grant N
N2014079 33 (2007-2009).

\end{document}